\documentclass[12pt]{article}

\usepackage{amsmath}

\usepackage{amsthm}

\usepackage{amsfonts}

\usepackage{amssymb}

\usepackage{amstext}

\usepackage{amsopn}

\usepackage{color}

\usepackage{ulem}

\usepackage{ulem}

\usepackage{geometry}
\def\C{\mathbb C}

\def\D{\mathbb D}

\def\limsup{\mathop{\overline{\rm lim}}}

\newtheorem{thm}{Theorem}

\newtheorem*{thmA}{Theorem~A}

\newtheorem{corollary}{Corollary}

\newtheorem{lemma}{Lemma}

\newtheorem{remark}{Remark}

\newtheorem{example}{Example}

\begin{document}

\title{\bf  Exact meromorphic solutions of Schwarzian differential equations}

\author{\textsc{Liangwen Liao and Chengfa Wu}}

\date{}

\maketitle

\begin{figure}[b]
\rule[-2.5truemm]{5cm}{0.1truemm}\\[2mm]
{\footnotesize  {\it Mathematics Subject Classification (2020)}:
Primary 34M05; Secondary  30D35.\par {\it Key words and phrases.} differential equation,  exact solution, Schwarzian
  differential equation, elliptic function. \par} 

\end{figure}

\begin{quotation}
\noindent{\sc\bf Abstract: This paper studies  exact meromorphic solutions of the autonomous Schwarzian differential equations. All transcendental meromorphic solutions of five canonical types (among six) of the autonomous Schwarzian differential equations are constructed explicitly. In particular, the solutions of four types are shown to be elliptic functions. Also, all transcendental meromorphic solutions that are locally injective or possess a Picard exceptional value are characterized for the remaining canonical type.}
\end{quotation}

\section{Introduction and Lemmas}

The Schwarzian derivative of a meromorphic function $f$ is
defined as
$$S(f,z) = \left({{f''}\over{f'}}\right)' - {1\over 2}\left({{f''}\over{f'}}\right)^2 = {{f'''}\over{f'}} -
{3\over 2}\left({{f''}\over{f'}}\right)^2.$$
It is well-known that $S(f,z) \equiv 0$ if and only if $f$ is a M\"{o}bius transformation. This property reveals that the Schwarzian derivative $S(f,z)$ measures how much $f$ differs from being a M\"{o}bius transformation. Another basic property of the Schwarzian derivative is that it is invariant under the M\"{o}bius group in the sense that $S(f,z)=S(\gamma \circ f,z)$, where $\gamma$ can be any M\"{o}bius transformation. The converse is also true, namely, if $S(g,z)=S(f,z)$, where $f,g$ are meromorphic functions, then there exits a  M\"{o}bius transformation $\gamma$ such that $g = \gamma \circ f$.

The Schwarzian derivative plays an essential role in various branches of complex analysis  \cite{Hille-book,Lehto1987,Steinmetz-1981} including univalent functions and conformal mappings.  It has also been shown that the Schwarzian derivative has close connections with second-order linear differential equations \cite{Laine-book} and Lax pairs of certain integrable partial differential equations \cite{Weiss1983}. In particular, it appears in the differential equation 
\begin{equation}\label{sch-de}
S(f,z)^p
= R(z,f) ={{P(z,f)}\over {Q(z,f)}},
\end{equation}
where $p$ is a positive integer, and $R(z,f)$ is an irreducible rational function in $f$ with meromorphic coefficients. The equation \eqref{sch-de} is known as the {\it Schwarzian
  differential equation}. Ishizaki \cite{Ishizaki-1991} obtained some Malmquist-type theorems of this equation and results concerning the deficiencies of its meromorphic solutions.  The growth of meromorphic solutions of the equation \eqref{sch-de} with polynomial coefficients has been studied by Liao and Ye \cite{Liao-Ye-2005}. A more complicated Schwarzian
type  differential equation was considered by Hotzel and Jank \cite{Hotzel1996Jank}.
If we restrict ourselves to the autonomous Schwarzian differential equation
 \begin{equation}\label{sch-de1}
S(f,z)^p
= R(f) ={{P(f)}\over {Q(f)}},
\end{equation}
where $P, Q$ are co-prime polynomials with constant coefficients,   Ishizaki \cite{Ishizaki-1991} obtained a Malmquist-Yosida-type result in which he gave a complete classification of the equation \eqref{sch-de1} possessing transcendental meromorphic solutions.

\begin{thmA} \label{autonomous Schwarzian DE}
 Suppose that the autonomous Schwarzian differential equation (\ref{sch-de1}) admits a transcendental
meromorphic solution. Then for some M\"{o}bius transformation $u = (af +b)/(cf +d), ad-bc\not=0,$
(\ref{sch-de1})  reduces into one of the following types
\begin{eqnarray}
S(u,z)&=&c\frac{(u-\sigma_1)(u-\sigma_2)(u-\sigma_3)(u-\sigma_4)}{(u-\tau_1)(u-\tau_2)(u-\tau_3)(u-\tau_4)} \label{ss-eq1}
\\
S(u,z)^3&=&c\frac{(u-\sigma_1)^3(u-\sigma_2)^3}{(u-\tau_1)^3(u-\tau_2)^2(u-\tau_3)} \label{ss-eq2}
\\
S(u,z)^3&=&c\frac{(u-\sigma_1)^3(u-\sigma_2)^3}{(u-\tau_1)^2(u-\tau_2)^2(u-\tau_3)^2} \label{ss-eq3}
\\
S(u,z)^2&=&c\frac{(u-\sigma_1)^2(u-\sigma_2)^2}{(u-\tau_1)^2(u-\tau_2)(u-\tau_3)} \label{ss-eq4}
\\
S(u,z)&=&c\frac{(u-\sigma_1)(u-\sigma_2)}{(u-\tau_1)(u-\tau_2)} \label{ss-eq5}
\\
S(u,z)&=&c \label{ss-eq6}
\end{eqnarray}
where $c \in \C, \tau_j $ are distinct constants, and $\sigma_j $ are constants, not necessarily distinct,  $j = 1,\dots,4$.
\end{thmA}

\begin{remark}
We remark that the conclusion of Theorem A does not hold for rational solutions of equation (\ref{sch-de1}). For instance, the function $$f(z) = -\frac{3 }{2 ( z+a)^2},$$ where $a$ is an arbitrary constant, satisfies the equation $S(u,z)=u$ but it cannot be transformed into any type of \eqref{ss-eq1}-\eqref{ss-eq6} via M\"{o}bius transformations. It is also noted that $f$ can be viewed as a fixed point of the Schwarzian operator and we refer the readers to the reference \cite{Zemyan2011} for the details on fixed points and $N$-cycles of the Schwarzian operator.
\end{remark}

The above theorem intimates  that to study the autonomous Schwarzian differential equation (\ref{sch-de1}), it suffices to consider the equations \eqref{ss-eq1}--\eqref{ss-eq6}. We will show that all transcendental meromorphic solutions of the equations \eqref{ss-eq1}-\eqref{ss-eq4} are elliptic functions and can be explicitly constructed. 
It is also shown that all transcendental meromorphic solutions of the equation (\ref{sch-de1}) can be characterized by imposing some conditions on them. The precise statements of these results are as follows.

\begin{thm} \label{Schwarzian DE-Picard exceptional vaule}
If the Schwarzian differential equation (\ref{sch-de1}) admits a transcendental meromorphic solution $f$ with a Picard exceptional value   $\xi\in\hat{\mathbb{C}}$, then by some M$\ddot{o}$bius transformation $f = \gamma_1(u)$,
(\ref{sch-de1})  reduces into either
$$S(u,z)=c\frac{(u-\sqrt2i)(u+\sqrt2i)}{(u-1)(u+1)},$$
and the transcendental meromorphic solutions of (\ref{sch-de1}) are $f(z)=\gamma_1(\sin(\alpha z+\beta))$, where $\alpha=  \sqrt{2c}  $ and $\beta$ is a constant; or
$$S(u,z)=c,$$
and  all   solutions of (\ref{sch-de1}) are
 $f(z)=\gamma_2 (e^{\alpha z})$, where $\alpha=  \sqrt{-2c}$ and $\gamma_2$ is any M\"{o}bius transformation.
\end{thm}

\begin{remark}
Theorem \ref{Schwarzian DE-Picard exceptional vaule} shows that any transcendental meromorphic solution of equations \eqref{ss-eq1}-\eqref{ss-eq4} must have infinitely many poles.
\end{remark}

The result below follows immediately from Theorem \ref{Schwarzian DE-Picard exceptional vaule}.

\begin{corollary} 
If the Schwarzian differential equation (\ref{sch-de1}) admits a transcendental entire solution $f$, then by some linear transformation $f = L_1(u)$,
(\ref{sch-de1})  reduces into either
$$S(u,z)=c\frac{(u-\sqrt2i)(u+\sqrt2i)}{(u-1)(u+1)},$$
and the entire solutions of (\ref{sch-de1}) are $f(z)=L_1(\sin(\alpha z+\beta))$, where $\alpha=  \sqrt{2c}  $ and $\beta$ is a constant; or
$$S(u,z)=c,$$
and  all entire solutions of (\ref{sch-de1}) are
 $f(z)=L_2(e^{\pm \alpha z})$, where $\alpha=  \sqrt{-2c}$ and $L_2$ is any linear transformation.
\end{corollary}

\begin{thm} \label{locally injective case}
If the Schwarzian differential equation (\ref{sch-de1}) admits a locally injective transcendental meromorphic solution, then by some M\"{o}bius transformation $f = \gamma(u)$,  \eqref{sch-de1}  reduces into
$$S(u,z)=c,$$
 and  all solutions of (\ref{sch-de1}) are
 $f(z)=\gamma(e^{\alpha z})$, where $\alpha=  \sqrt{-2c}$ and $\gamma$ is any M\"{o}bius transformation.
\end{thm}

Rewrite the equation \eqref{ss-eq1} as
\begin{eqnarray}
S(u,z)&=&c\frac{(u-\sigma_1)(u-\sigma_2)(u-\sigma_3)(u-\sigma_4)}{(u-\tau_1)(u-\tau_2)(u-\tau_3)(u-\tau_4)} \nonumber
 \\
 &=& \frac{r_4 u^4 + r_3 u^3 + r_2 u^2 + r_1 u + r_0 u }{(u-\tau_1)(u-\tau_2)(u-\tau_3)(u-\tau_4)}, \label{ss-eq1-new form}
\end{eqnarray}
and denote by
 \begin{eqnarray} \label{e_i}
 e_1 = \sum_{j=1}^4 \tau_j, \quad e_2 = \sum_{1\leq j < k \leq 4} \tau_j \tau_k, \quad e_3 =   \sum_{1\leq j < k < l \leq 4} \tau_j \tau_k \tau_l, \quad e_4 = \prod_{j=1}^4 \tau_j.
 \end{eqnarray}
 Then we can construct all transcendental meromorphic solutions to the  equation \eqref{ss-eq1}.

\begin{thm} \label{classification of schwarzian DE1}
All transcendental meromorphic solutions of the  equation  \eqref{ss-eq1-new form}  are elliptic functions   of the form
\begin{eqnarray}\label{solution to ss-eq1}
 u(z) 
 &=& a-\frac{b}{\wp (z-z_0;g_2,g_3)-d},
\end{eqnarray}
where $\wp(z;g_2,g_3)$ is the Weierstrass elliptic function,  $z_0\in\C$ is arbitrary, $a=\tau_i $ and $b,d, g_2,g_3$ are constants that depend on $c$, $\sigma_i$ and $\tau_i,i=1,2,3,4$. Further, with $e_i \, (i=1,2,3,4)$ defined in \eqref{e_i} and
$$q_i =    \prod_{\substack{1\leq j   \leq 4 \\ j\not = i}}   (\tau_i - \tau_j), \quad i=1,2,3,4,$$
the equation  \eqref{ss-eq1-new form}  admits solutions of the form \eqref{solution to ss-eq1} if and only if the following parameter relations hold
  \begin{eqnarray}
 r_0 &=& \frac{b}{2q_i} \left(3 e_3^2-8 e_2 e_4\right), \quad r_1 = \frac{2b}{q_i} \left(6 e_1 e_4-e_2 e_3\right), \quad  r_2  = \frac{b}{q_i} \left(2 e_2^2-3 e_1 e_3-24 e_4\right), \label{parameter constraint-begin} \\
 r_3 &=& \frac{2b}{q_i} \left(6 e_3-e_1 e_2\right), \quad  r_4 = \frac{b}{2q_i} \left(3e_1^2- 8e_2\right), \label{parameter constraint-2}\\ 
 d &=& \frac{b}{6q_i}   \left[\sum_{\substack{1\leq j < k \leq 4 \\ j,k \not = i}} (\tau_j - \tau_k)^2 - \sum_{\substack{1\leq j   \leq 4 \\ j\not = i}} 2 (\tau_i - \tau_j)^2\right], \\
 g_2 &=& \frac{4b^2}{3q_i^2} \left(e_2^2-3 e_1 e_3+12 e_4\right), \quad
 g_3 = \frac{4b^3}{27q_i^3} \left(2 e_2^3-9 e_1 e_2 e_3-72 e_2 e_4+27 e_3^2+27 e_1^2 e_4\right), \label{parameter constraint-end}\\
 \Delta &=& g_2^3-27g_3^2 = \frac{16 b^6}{q_i^6} \prod_{1\leq j < k \leq 4} (\tau_j - \tau_k)^2  \not=0.
\end{eqnarray}

\end{thm}

\begin{remark}
Theorem \ref{classification of schwarzian DE1} indicates that  the equation \eqref{ss-eq1-new form} has transcendental meromorphic solutions only if the parameters $c, \sigma_i, \tau_i, i=1,2,3,4$ satisfy the conditions \eqref{parameter constraint-begin} and \eqref{parameter constraint-2}. In addition, the solution \eqref{solution to ss-eq1} has just one free parameter, which implies that the general solution of equation \eqref{ss-eq1-new form}  should have more complicated singularities other than poles.
\end{remark}

In view of the invariance of Schwarzian derivatives under the M\"{o}bius group, we may compose the solution $u$ of  equations \eqref{ss-eq2}-\eqref{ss-eq4} with a M\"{o}bius transformation such that    $ \tau_1, \tau_2 $ and $ \tau_3 $ can be any distinct desired numbers, and this allows us to derive all transcendental meromorphic solutions to the  equation equations \eqref{ss-eq2}-\eqref{ss-eq4}  explicitly.

\begin{thm} \label{classification of schwarzian DE2}
Let $ \tau_1 = 4, \tau_2=-3, \tau_3 =0$, then all transcendental meromorphic solutions to the  equation \eqref{ss-eq2} are elliptic functions.
Moreover, these solutions exist if and only if
\begin{eqnarray*}
 && \{\sigma_1, \sigma_2\} = \left\{  \sqrt{5} i , -  \sqrt{5} i \right\}.  
\end{eqnarray*}
and in this case, all the  transcendental meromorphic solutions to the  equation \eqref{ss-eq4} are given by
\begin{equation*}
u(z) =  -\frac{3 c}{c-74088 \wp \left(z-z_0;g_2,g_3\right)^3},
\end{equation*}
where  $\wp(z;g_2,g_3)$ is the Weierstrass elliptic function with $g_2 = 0,    g_3 = c/10584$, and $z_0\in\C$ is arbitrary.

\end{thm}

\begin{thm} \label{classification of schwarzian DE3}
Let $\{ \tau_1, \tau_2, \tau_3 \} = \{0,1,-1\}$, then all transcendental meromorphic solutions to the  equation \eqref{ss-eq3} are elliptic functions.
Moreover, these solutions exist if and only if  \begin{eqnarray*}
  \{\sigma_1, \sigma_2\} = \left\{\frac{i}{\sqrt{3}},-\frac{i}{\sqrt{3}}\right\},
\end{eqnarray*}
and in this case, all the  transcendental meromorphic solutions to the  equation \eqref{ss-eq3} are given by
\begin{equation*}
u(z) =  \frac{9 \left[9 \wp \left(z-z_0;g_2,g_3\right)+L^2\right] \wp'\left(z-z_0;g_2,g_3\right)}{2 L \left[81 \wp \left(z-z_0;g_2,g_3\right)^2-9L^2\wp \left(z-z_0;g_2,g_3\right)+L^4\right]}
\end{equation*}
where $L^{6} = -27c/64$, $\wp(z;g_2,g_3)$ is the Weierstrass elliptic function with $g_2 = 0,  g_3 = c/432$, and $z_0\in\C$ is arbitrary.

\end{thm}

\begin{thm} \label{classification of schwarzian DE4}
Let $ \tau_1 = 0, \tau_2=1, \tau_3 =-1$, then all transcendental meromorphic solutions to the  equation \eqref{ss-eq4} are elliptic functions.
Moreover, these solutions exist if and only if  \begin{eqnarray*}
  \{\sigma_1, \sigma_2\} = \left\{\frac{i}{2},-\frac{i}{2}\right\},
\end{eqnarray*}
and in this case, all the  transcendental meromorphic solutions to the  equation \eqref{ss-eq4} are given by
\begin{equation*}
u(z) =  -\frac{1}{2 L}\frac{\left(8 \wp \left(z-z_0;g_2,g_3\right)+L^2\right)^2 \wp '\left(z-z_0;g_2,g_3\right)}{ \wp \left(z-z_0;g_2,g_3\right) \left(64 \wp \left(z-z_0;g_2,g_3\right)^2+L^4\right)},
\end{equation*}
where $c=   9 L^4/4$, $\wp(z;g_2,g_3)$ is the Weierstrass elliptic function with $g_2 = - c/36, g_3 = 0$, and $z_0\in\C$ is arbitrary.

\end{thm}

\begin{remark}
 It follows from Theorems \ref{Schwarzian DE-Picard exceptional vaule}-\ref{classification of schwarzian DE4} that all transcendental meromorphic solutions of the canonical Schwarzian differential equations \eqref{ss-eq1}-\eqref{ss-eq6} have been derived, except the solutions of equation \eqref{ss-eq5} that have no Picard exceptional values. Although we are not able to prove that any transcendental meromorphic solution of equation \eqref{ss-eq5} must have Picard exceptional value(s), we conjecture this is true.
\end{remark}

\section{Preliminaries}

The important tools in our proof include Wiman-Valiron theorem and Wiman-Valiron theory. 
Let  $f$ be  a transcendental entire function, and write
$$f(z) = \sum_{n=0}^{\infty}a_nz^n.$$
As usual, for $r>0$, we denote the maximum term by $\mu(r,f)$, the central index by $\nu(r,f)$, and the maximum modulus by $M(r,f)$, i.e.,
$$\mu(r,f) = \max_{|z|=r}|a_nz^n|, \ \  \nu(r,f) = \sup\{n| |a_n|r^n=\mu(r,f)\},
\ \ M(r,f) = \max_{|z|=r}|f(z)|.$$
\begin{lemma}[Wiman-Valiron Theorem \cite{Bergweiler-2001}] \label{WV-thm}
 There exists a set   $F\subset[1,+\infty)$ satisfying
$$\int_F\frac{d t}{t}<\infty$$
with the following property: if $(z_k)$ is a sequence in $\mathbb{C}$ with $|f(z_k)|=M(|z_k|,f), |z_k|\not\in F$ and
$z_k\to\infty$, and if $\nu_k=\nu(|z_k|,f),$ then
$$\frac{f\left(z_k+\dfrac{z_k}{\nu_k}z\right)}{f(z_k)}\to e^z$$
as $k\to\infty.$
\end{lemma}

\begin{lemma}[\cite{Laine-book}] \label{Wiman-Varilon-lemma}

Let $f$ be a transcendental entire
function,   $0<\delta<{1\over 4}$ and $|z|=r$ such that
$$|f(z)|>M(r,f){\nu(r,f)}^{-{1\over 4}+\delta}\eqno(2.2)$$
holds. Then there exists a set $F \subset (0, + \infty)$ of finite logarithmic measure, i.e., $\int_F dt/ t< +\infty$ such that
$$f^{(m)}(z) = {\left({{\nu(r,f)}\over z}\right)}^m(1+o(1))f(z)\eqno(2.3)$$
holds for all $m\ge 0$ and all $r\not\in F$.
\end{lemma}

The Schwarzian derivative has a fundamental relation with second-order linear ordinary differential equations.

\begin{lemma}{\cite[p.~110]{Laine-book}} \label{Schwarzian-2nd order linear ODE}
 Let $A(z)$ be analytic in a simply connected domain $\Omega$. Then, for any two linearly independent solutions $f_{1}, f_{2}$ of
 \begin{equation}\label{2nd order linear ODE}
    f''(z)+A(z)f(z) = 0,
 \end{equation}
    their quotient $g=f_{1} / f_{2}$ is locally
injective and satisfies the differential equation
\begin{equation}\label{ratio}
  S(g,z)=2 A(z).
\end{equation}
Conversely, let $g$ be a locally injective meromorphic function in $\Omega$ and define $A(z)$
by \eqref{ratio}. Then $A(z)$ is analytic in $\Omega$ and the differential equation \eqref{2nd order linear ODE} admits two linearly independent solutions $f_{1}, f_{2}$ such that $g=f_{1} / f_{2}$.
\end{lemma}

\begin{remark}
The  lemma above has crucial applications in differential equations. In particular, it has been used by Bergweiler and Eremenko \cite{Bergweiler2017Eremenko} to solve the Bank-Laine conjecture, which concerns the zero distribution of solutions of   equation \eqref{2nd order linear ODE} where $A$ is an entire function of finite order.
\end{remark}

Now we introduce some terminologies in Nevanlinna theory \cite{Laine-book}. Let $f$ be a meromorphic function on $\C$ and $n(r,f)$ denote the number of poles of $f$ in the disk $\D(r) = \{z\in \C ||z| <r\}$, counting multiplicity. The Nevanlinna characteristic function of $f$ is defined as
\begin{equation*}
  T (r, f ) = m(r, f ) + N(r, f ),
\end{equation*}
where
\begin{eqnarray*}
  m(r, f) &=& \int_{0}^{2 \pi} \log ^{+} |f\left(r e^{i \theta}\right)| \frac{d \theta}{2 \pi}, \\
  N(r, f)&=&  n(0, f) \log r+\int_{0}^{r}\left[n(t, f)-n(0, f)\right] \frac{d t}{t},
\end{eqnarray*}
with $\log ^{+} x=\max \{0, \log x\}$. We note that $m(r, f)$ and $N(r, f)$ are called the proximity function and integrated counting function, respectively.
Next, we  define the order of  $f$  by
\begin{equation*}
  \rho(f)=\limsup _{r \rightarrow \infty} \frac{\log T(r, f)}{\log r}.
\end{equation*}
The following result of Liao and Ye \cite{Liao-Ye-2005} says that   the order of meromorphic solutions of equation \eqref{sch-de1} is bounded from above by $2$.

\begin{lemma}  \label{growth}
Let $f$ be a meromorphic solution of the autonomous Schwarzian differential equation \eqref{sch-de1},
then $\rho (f ) \leq 2$.
\end{lemma}

\section{Proof of main results}

We first recall the definition of totally ramified values: we call a point $a \in \overline{\C}$  a totally ramified value of a meromorphic function $f$ if all $a$-points of $f$ are multiple. According to a classical result of Nevanlinna, a non-constant function meromorphic in the plane can have at most four totally ramified values while a non-constant entire function can have at most two finite totally ramified values. We also need the following results.

\begin{lemma}[\cite{Laine-book}]\label{log-de-lemma}
Let $f(z)$ be
a nonconstant meromorphic function. Then $$m\left(r,\frac{f'}{f}\right)=O(\log
r),   $$if $f$ is of finite order, and
$$m\left(r,\frac{f'}{f}\right)=O(\log (r T(r,f))),$$  possibly
outside a set $E$ of $r$ with finite linear measure, if $f(z)$ is of
infinite order.
\end{lemma}
\begin{lemma}[\cite{Zhangx-liao}] \label{nonlinear-eq-thm1} If the differential equation
\begin{equation}\label{1stad-eq1}
w^2+R(z)(w^{\prime})^2= Q(z),
\end{equation}
where $R, Q$ are nonzero rational functions,  admits a transcendental meromorphic solution
$f$,
 then $Q \equiv A$ is a constant, the multiplicity of zeros of $R(z)$ is no greater than 2 and $f(z)=
 \sqrt{A}\cos\alpha(z)$, where  $\alpha(z)$ is a primitive of
  $1/\sqrt{R(z)}$
    such that
$
 \sqrt{A}\cos\alpha(z)$ is a transcendental
meromorphic function.
\end{lemma}

\subsection{Proof of Theorem \ref{Schwarzian DE-Picard exceptional vaule}}

Let $f$ be a transcendental meromorphic solution with a Picard exceptional value of the equation (\ref{sch-de1}). It follows from Theorem A that by some M\"{o}bius transformation
$$u = \dfrac{af +b}{cf +d}, \quad ad-bc\not=0,$$
$u$ satisfies one of the equations \eqref{ss-eq1}-\eqref{ss-eq6}.

If $u$ satisfies  the equation \eqref{ss-eq1}, then $u$ has four   totally ramified values $\tau_1,\tau_2,\tau_3,\tau_4$. This is impossible since $u$ has a Picard exceptional value. If $u$ satisfies  the equation (\ref{ss-eq2}), then $u$ has three totally  ramified values $\tau_1,\tau_2,\tau_3$. Thus,   the Picard exceptional value of $u$ must be one of  them.
Without loss of generality, we may assume $\tau_3$ is  a Picard exceptional value of $u$. Let
$$v=\dfrac{1}{u-\tau_3},$$
then $v$ has at most finitely many poles   and satisfies the following differential equation
\begin{equation}\label{ss-eq-2-1}
S(v,z)=c'\frac{(v-\sigma'_1)^3(v-\sigma'_2)^3}{(v-\tau'_1)^3(v-\tau'_2)^2}
\end{equation}
Assume $\zeta_1,\cdots,\zeta_n$ are the poles (counting multiplicities) of $v$, then
$v(z)=g(z)/P(z)$, where $g(z)$ is a transcendental entire function and $P(z)=(z-\zeta_1)\cdots(z-\zeta_n).$
We choose $z_k\to\infty$ such that $|z_k|\not\in F$ and $|g(z_k)|=M(|z_k|,g).$ Let
$$h_k(z)=\frac{v(z_k+\rho_k z)}{v(z_k)},$$
 where $\displaystyle \rho_k=\frac{z_k}{\nu_k}, \nu_k=\nu(|z_k|,g)$, then by Lemma \ref{WV-thm}, we have
$$\lim_{k\to\infty}h_k(z)=\lim_{k\to\infty}\frac{v(z_k+\rho_k z)}{v(z_k)}
=\lim_{k\to\infty}\frac{g(z_k+\rho_kz)}{g(z_k)}\frac{P(z_k)}{P(z_k+\rho_kz)}=e^z.$$
Thus
$$\lim_{k\to\infty}\frac{\rho_k v'(z_k+\rho_kz)}{v(z_k)}=\lim_{k\to\infty}h_k'(z)=e^z,$$
and
$$\lim_{k\to\infty}\frac{\rho_k^2v''(z_k+\rho_kz)}{v(z_k)}=\lim_{k\to\infty}h_k''(z)=e^z.$$
It follows from (\ref{ss-eq-2-1}) that
\begin{equation}\label{scale1}
\frac1{v(z_k)}\left(\frac{1}{\rho_k}\right)^2\left(\frac{h_k'''(z)}{h_k'(z)}-\frac32
\left(\frac{h_k''(z)}{h_k'(z)}\right)^2\right)=    c'
\frac{(h_k(z)-\sigma_1'/v(z_k))^3(h_k(z)-\sigma_2'/v(z_k))^3}
{(h_k(z)-\tau_1'/v(z_k))^3(h_k(z)-\tau_2'/v(z_k))^2}.
\end{equation}
Noting the selection of $z_k$, we have $$\displaystyle \lim_{k\to\infty}\frac{\nu_k^M}{v(z_k)}=0$$
for any positive number $M$. Thus, the left side of the  equation \eqref{scale1} tends  to zero while the right side of equation \eqref{scale1} tends tends to $c'e^z$ as $k\to\infty$, which is a contradiction. Thus $u$ cannot satisfy (\ref{ss-eq2}). With similar arguments, we can prove that  $u$ satisfies neither (\ref{ss-eq3}) nor (\ref{ss-eq4}).

If $u$ satisfies the equation (\ref{ss-eq5}), then $u$ has two    totally ramified values $\tau_1,\tau_2$. Then we distinguish two cases.\\
 Case 1: one of $\tau_1$ and $\tau_2$ is the Picard exceptional value of $u$, by the same arguments as above, we get a contradiction. \\
 Case 2: both of $\tau_1$ and $\tau_2$ are not the Picard exceptional value of $u$. Without loss of generality, we may assume the Picard exceptional value of $u$ is infinity.  Otherwise,
we may consider a composition of a M\"{o}bius transformation and the function $u$. Thus we
can express $u$ as
$$u(z) =\frac{g(z)}{P(z)},$$
where $g(z)$ is a transcendental entire function and $P(z)$ is a polynomial. For any $r >0,$ let
$$|g(z_0)| = M(g, r), \quad |z_0| = r.$$ Then, by Lemma \ref{Wiman-Varilon-lemma}, there exists a set $F \subseteq (0, +\infty)$ with a finite logarithmic
measure such that
$$\frac{u'(z_0)}{u(z_0)}=\frac{g'(z_0)}{g(z_0)}-\frac{P'(z_0)}{P(z_0)}
= \frac{\nu(g, r)}{z_0}(1+o(1)),$$
$$\frac{u''(z_0)}{u(z_0)}=\frac{g''(z_0)}{g(z_0)}-\frac{P''(z_0)}{P(z_0)}-
2\frac{u'(z_0)}{u(z_0)}\frac{P'(z_0)}{P(z_0)}
= \left(\frac{\nu(g, r)}{z_0}\right)^2(1+o(1)),$$
and
$$\frac{u'''(z_0)}{u(z_0)}=\frac{g'''(z_0)}{g(z_0)}-\frac{P'''(z_0)}{P(z_0)}-
3\frac{u''(z_0)}{u(z_0)}\frac{P'(z_0)}{P(z_0)}-3\frac{u'(z_0)}{u(z_0)}\frac{P''(z_0)}{P(z_0)}
= \left(\frac{\nu(g, r)}{z_0}\right)^3(1+o(1)),$$
for all sufficiently large $r\not\in F.$ Thus the equation (\ref{ss-eq5})  becomes
$$\left(\frac{\nu(g, r)}{z_0}\right)^2(1+o(1))-\frac32\left(\frac{\nu(g, r)}{z_0}(1+o(1))\right)^2=c'(1+o(1)).$$
This leads to $$\nu(r,g)\sim Ar \text{ and }\rho(g)=1.$$ Hence $\rho(u)=1.$

 By computing the Laurent expansions on both sides of  (\ref{ss-eq5}), we may obtain

\begin{itemize}
  \item  $u'(z)=0$ if and only if $u(z)=\tau_1$ or $u(z)=\tau_2$.
  \item all the zeros of $u'$ are simple.
\end{itemize}
Without loss of generality, we may assume $\tau_1=1,\tau_2=-1. $ Thus, $$\displaystyle \frac{(u')^2}{u^2-1}$$ is a meromorphic function having only finitely many poles and no zeros. Noting $\rho(u)=1$, we have
$$\displaystyle Q^2(z) \frac{(u')^2}{u^2-1}=e^{h(z)},$$
where $Q(z)$ is a nonzero polynomial with simple zeros and $h(z)$ is an entire function. Then, by Lemma \ref{log-de-lemma}, we have
\begin{eqnarray*}
T(r, e^{h})&=&m(r, e^{h})
\\
&\leq& 2 m\left(r,Q\right)  +m\left(r,\frac{u'}{u-1}\right)+m\left(r,\frac{u'}{u+1}\right)
\\
&=&O(\log r).
\end{eqnarray*}
This implies $e^h$ is a polynomial and hence $h$ is a constant. Without loss of generality, we may assume $h = 1$, then $u$ satisfies the differential equation
\begin{equation}\label{sin-eq}
u^2 -Q(z)^2(u')^2=1.
\end{equation}
If $\deg Q\ge 1$, then by the equation (\ref{sin-eq}) and Lemma \ref{Wiman-Varilon-lemma}, we have $$\nu(r,u)\sim Ar^{1-\frac{\deg P}2},$$ where $A$ is a positive number, but this contradicts with $\rho(u)=1.$ Hence $Q(z)$ is a constant. It is easy to see that the solutions of (\ref{sin-eq}) are of the form
$$u=\sin(\alpha z+\beta),$$
where $\alpha, \beta$ are constants with $Q^2\alpha^2=-1.$ Substituting $u=\sin(\alpha z+\beta)$ into (\ref{ss-eq5}) and noting $\tau_1=1, \tau_2=-1$, we obtain that $$\alpha=   \sqrt{2c}, \quad \sigma_1=\sqrt2i, \quad \sigma_2=-\sqrt2i.$$ Thus we get the conclusion.

Finally, if $u$ satisfies  equation (\ref{ss-eq6}), then $R(f)$ must be a constant, say A, and hence $c^p = A$.
It is easy to check that $u(z)=e^{\alpha z}$ is a solution of the equation (\ref{ss-eq6}), where $\alpha=\sqrt{-2c}$. Then it follows  from the invariance property  of the Schwarzian derivative under M\"{o}bius transformations that  all the solutions of (\ref{ss-eq6}) are given by $u=\gamma(e^{\alpha z})$, where $\gamma$ is a M\"{o}bius transformation and $\alpha=\sqrt{-2c}.$ Hence, in this case,  all the solutions of the equation (\ref{sch-de1})
are $f(z)=\gamma(e^{\alpha z})$, where $\gamma$ is a M\"{o}bius transformation and $\alpha=\sqrt{-2}A^{\frac1{2p}}.$

\subsection{Proof of Theorem  \ref{locally injective case}}

Suppose $f$ is a locally injective transcendental meromorphic solution  of the equation \eqref{sch-de1}. According to Theorem A, there exits a M\"{o}bius transformation $\gamma_1$ such that $u=\gamma_1(f)$ is also a locally injective transcendental meromorphic function and satisfies one of the equations \eqref{ss-eq1}-\eqref{ss-eq6}. Then it follows from Lemma \ref{Schwarzian-2nd order linear ODE} that $S(u,z)$ is entire. This  implies   $u$ cannot satisfy any of the equations \eqref{ss-eq1}-\eqref{ss-eq5}. Otherwise, $u$ has at least one Picard exceptional value. By Theorem  \ref{Schwarzian DE-Picard exceptional vaule}, it indicates that $u = \gamma_2(\sin(\alpha z+\beta))$, where $\alpha, \beta$ are constants and $\gamma_2$ is a M\"{o}bius transformation. Nevertheless, this contradicts with the fact  that $u$ is locally injective. As a consequence, $u$ can only satisfy equation \eqref{ss-eq6} and then the conclusion follows immediately from Theorem \ref{Schwarzian DE-Picard exceptional vaule}.

\subsection{Proof of Theorem  \ref{classification of schwarzian DE1}}

Suppose $u$ is a transcendental meromorphic solution to the  equation \eqref{ss-eq1-new form}, then Theorem \ref{Schwarzian DE-Picard exceptional vaule} shows that $u$ must have infinitely many poles. By comparing the Laurent expansions on both sides of  the  equation \eqref{ss-eq1-new form}, we deduce that all the poles of  $u$  are simple and all the poles (if they exist) of $S(u,z)$ come from the zeros of $u'$. Since all the poles of $S(u,z)$ are double, it follows that all  zeros of $u'$ should be simple, and at any zero of $u'$, $u(z)$ assumes one of the $\tau_i, i=1,2,3,4$. This means any zero of $u-\tau_i$ must be double. Therefore,
\begin{equation}
G(z) = \frac{u'^2}{(u-\tau_1)(u-\tau_2)(u-\tau_3)(u-\tau_4)}
\end{equation}
is a nonvanishing entire function, and there exists an entire function $g(z)$ such that $G=e^{g}$.  According to Theorem \ref{growth}, $u$ has finite order of growth.
Then
 we have
\begin{eqnarray*}
T(r, e^{g})&=&m(r, e^{g})
\\
&\leq&  m\left(r,\frac{u'}{(u-\tau_1)(u-\tau_2)}\right)+m\left(r,\frac{u'}{(u-\tau_3)(u-\tau_4)}\right)
\\
&\leq& \sum_{i=1}^4 m\left(r,\frac{u'}{u-\tau_i}\right) + O(1)
\\
&=&O(\log r),
\end{eqnarray*}
where the last equality follows from Lemma \ref{log-de-lemma}. This implies $e^g$ is a polynomial and hence $g = C$ is a constant. As a consequence, $u$ satisfies the differential equation
\begin{equation*}
u'^2 = K (u-\tau_1)(u-\tau_2)(u-\tau_3)(u-\tau_4), \quad  K = e^C
\end{equation*}
whose general solution is given by \cite{Conte2015NgWu}
\begin{eqnarray}
 u(z) &=& K^{-1 / 2}\left(A-\frac{\wp'(w;g_2,g_3)}{\wp(z-z_0;g_2,g_3)-\wp(w;g_2,g_3)}\right) \nonumber
 \\
 &=& a-\frac{b}{\wp (z-z_0;g_2,g_3)-d} \label{solution to Schwarzian DE1}
\end{eqnarray}
where $\wp(z;g_2,g_3)$ is the Weierstrass elliptic function, $z_0\in\C$ is arbitrary and $a,b,d,g_2,g_3$ are constants that depend on $K$ and $\tau_i,i=1,2,3,4$. Finally, by substituting \eqref{solution to Schwarzian DE1} into \eqref{ss-eq1-new form} and applying the differential equation satisfied by $\wp(z;g_2,g_3)$
\begin{equation*}
  \wp'^2 = 4\wp^3 - g_2 \wp - g_3,
\end{equation*}
where $ \Delta = g_2^3 - 27 g_3^2 \not= 0$, it can be  computed that $a $ should be equal to one of the $ \tau_i, i=1,2,3,4$, and other parameters should satisfy the relations \eqref{parameter constraint-begin}-\eqref{parameter constraint-end}.

\subsection{Proof of Theorem  \ref{classification of schwarzian DE2}}

Let $u$ be a transcendental meromorphic solution to the  equation \eqref{ss-eq2}, then Theorem \ref{Schwarzian DE-Picard exceptional vaule} implies that $u$ must have infinitely many poles. With similar arguments as in Theorem  \ref{classification of schwarzian DE1}, we find that
\begin{itemize}
  \item  all the poles of  $u$  are simple;
  \item $u'(z) = 0$ if and only if $u(z) = \tau_i $  for some $i \in \{1,2,3\}$;
  \item if $u(z) = \tau_1 $, then $z$ is a simple zero of $u'$;
  \item if $u(z) = \tau_2 $, then $z$ is a double zero of $u'$;
  \item if $u(z) = \tau_3 $, then $z$ is a  zero of $u'$ of order $5$.
\end{itemize}
It follows that
\begin{equation}
G(z) = \frac{u'^6}{(u-\tau_1)^3(u-\tau_2)^4(u-\tau_3)^5}
\end{equation}
is a nonvanishing entire function, and hence, there exists an entire function $g(z)$ such that $G=e^{g}$.  According to Theorem \ref{growth}, $u$ has finite order of growth.
Then we have
\begin{eqnarray*}
T(r, e^{g})&=&m(r, e^{g})
\\
&\leq& m\left(r,\frac{u'^3}{(u-\tau_1)^3(u-\tau_2)^3(u-\tau_3)^3}\right)+ m\left(r,\frac{u'}{u-\tau_2}\right)+  m\left(r,\frac{u'^2}{(u-\tau_3)^2}\right)
\\
&\leq& 3m\left(r,\frac{u'}{u-\tau_1}\right)+4m\left(r,\frac{u'}{u-\tau_2}\right)+5m\left(r,\frac{u'}{u-\tau_3}\right)+ O(1)
\\
&=&O(\log r).
\end{eqnarray*}
This indicates that $g= C$ is a constant and hence $u$  satisfies the differential equation
\begin{equation}  \label{subequation-2}
u'^6 = K (u-\tau_1)^3(u-\tau_2)^4(u-\tau_3)^5, \quad  K =  e^C.
\end{equation}
Since the  elliptic curve parametrized by  $u$ and $u'$ has genus one, the general solution of the above equation should be elliptic functions. Let
\begin{equation}\label{subequation-3-u to v}
  u(z) = \frac{1}{v(z)} + \tau_3,
\end{equation}
then the equation \eqref{subequation-2} reduces to
\begin{equation}  \label{subequation-3-v}
v'^6 = K [(\tau_1-\tau_3)v - 1]^3[(\tau_2-\tau_3)v - 1]^4. 
\end{equation}
By using the singularity methods (see \cite{Conte2018NgWu,Ng2019Wu} and the references therein), we find that the general solution to \eqref{subequation-3-v} reads
\begin{equation}\label{subequation-3-solution of v}
 v(z) = h-\frac{23328 \left[6 \wp (z-z_0;g_2,g_3)^3+\wp '(z-z_0;g_2,g_3)^2\right]}{5 K (\tau_1-\tau_3)^3 (\tau_2-\tau_3)^4},
\end{equation}
where   $\wp(z;g_2,g_3)$ is the Weierstrass elliptic function with $g_2=0$,  $z_0\in\C$ is arbitrary and  $h,g_3$ are constants  depending on $K$ and $\tau_i,i=1,2,3$. Finally, with $ \tau_1=4, \tau_2=-3, \tau_3 =0 $, substituting \eqref{subequation-3-u to v} and \eqref{subequation-3-solution of v} into \eqref{ss-eq2} yields
 the solution of equation \eqref{ss-eq2}
\begin{equation*}
u(z) =  -\frac{3 c}{c-74088 \wp \left(z-z_0;0,g_3\right)^3},
\end{equation*}
where $  g_3 = c/10584$ and $z_0\in\C$ is arbitrary, provided that
\begin{eqnarray*}
 && \{\sigma_1, \sigma_2\} = \left\{  \sqrt{5} i , -  \sqrt{5} i \right\}.  
\end{eqnarray*}
This completes the proof.

\subsection{Proof of Theorem  \ref{classification of schwarzian DE3}}

Suppose $u$ is a transcendental meromorphic solution to the  equation \eqref{ss-eq3}, then Theorem \ref{Schwarzian DE-Picard exceptional vaule} implies that $u$ must have infinitely many poles. Using similar arguments as in Theorem  \ref{classification of schwarzian DE1}, we can show that
\begin{itemize}
  \item  all the poles of  $u$  are simple;
  \item $u'(z) = 0$ if and only if $u(z) = \tau_i $  for some $i \in \{1,2,3\}$;
  \item all the zeros of $u'$ are double.
\end{itemize}
It follows that
\begin{equation}
G(z) = \frac{u'^3}{(u-\tau_1)^2(u-\tau_2)^2(u-\tau_3)^2}
\end{equation}
is a nonvanishing entire function, and hence, there exists an entire function $g(z)$ such that $G=e^{g}$.  Since the order of $u$ is finite, we have
\begin{eqnarray*}
T(r, e^{g})&=&m(r, e^{g})
\\
&\leq& m\left(r,\frac{u'}{u-\tau_1}\right)+  m\left(r,\frac{u'}{(u-\tau_2)(u-\tau_3)}\right)+m\left(r,\frac{u'}{(u-\tau_1)(u-\tau_2)(u-\tau_3)}\right)
\\
&\leq& 2 \left[ \sum_{i=1}^3 m\left(r,\frac{u'}{u-\tau_i}\right) \right] + O(1)
\\
&=&O(\log r).
\end{eqnarray*}
This implies that $g= C$ is a constant and hence $u$  satisfies the differential equation
\begin{equation}  \label{subequation-3}
u'^3 = K (u-\tau_1)^2(u-\tau_2)^2(u-\tau_3)^2, \quad  K =  e^C.
\end{equation}
Since the  elliptic curve parametrized by  $u$ and $u'$ has genus one, the general solution of the above equation should be elliptic functions. By using the singularity methods, we find that the general solution to \eqref{subequation-3} can be expressed as
\begin{eqnarray} \label{solution to Schwarzian DE3}
 u(z) &=&  \frac{1}{L}     \left(\frac{\left(1+i \sqrt{3}\right) \left(\wp '\left(z-z_0;g_2,g_3\right)-A_1\right)}{4 \left(\wp \left(z-z_0;g_2,g_3\right)-B_1\right)}+\frac{\left(1-i \sqrt{3}\right) \left(\wp '\left(z-z_0;g_2,g_3\right)-A_2\right)}{4 \left(\wp \left(z-z_0;g_2,g_3\right)-B_2\right)}\right) \nonumber
 \\
&& + \frac{1}{3} (\tau_1+\tau_2+\tau_3)
\end{eqnarray}
where $L^3 = K$, $\wp(z;g_2,g_3)$ is the Weierstrass elliptic function,  $z_0\in\C$ is arbitrary and  $A_1,A_2,B_1,B_2,g_2,g_3$ are constants  depending on $K$ and $\tau_i,i=1,2,3$. Finally, with $\{ \tau_1, \tau_2, \tau_3 \} = \{0,1,-1\}$, substituting \eqref{solution to Schwarzian DE3} into \eqref{ss-eq3} yields that
\begin{eqnarray*}
  \{\sigma_1, \sigma_2\} = \left\{\frac{i}{\sqrt{3}},-\frac{i}{\sqrt{3}}\right\}, \quad A_1 = A_2 = g_2 = 0, \quad g_3 = \frac{c}{432},
  \\
  B_1 = \frac{1}{18} \left(1-i \sqrt{3}\right) L^2, \quad B_2 = \frac{1}{18} \left(1+i \sqrt{3}\right) L^2, \quad  L^{6} = - \frac{27}{64} c.
\end{eqnarray*}
In this case, the equation \eqref{ss-eq3} reduces to
\begin{eqnarray*}
   S(u,z)^3 = c \frac{(u^2 + 1/3)^3}{u^2 (u^2-1)^2},
\end{eqnarray*}
and     the solution \eqref{solution to Schwarzian DE3} becomes
\begin{equation*}
u(z) =  \frac{9 \left[9 \wp \left(z;g_2,g_3\right)+L^2\right] \wp'\left(z;g_2,g_3\right)}{2 L \left[81 \wp \left(z;g_2,g_3\right)^2-9L^2\wp \left(z;g_2,g_3\right)+L^4\right]}
\end{equation*}
where $L^{6} = -27c/64, g_2 = 0, g_3 = c/432$.

\subsection{Proof of Theorem  \ref{classification of schwarzian DE4}}

Let $u$ be a transcendental meromorphic solution to the  equation \eqref{ss-eq4}, then Theorem \ref{Schwarzian DE-Picard exceptional vaule} implies that $u$ must have infinitely many poles. With similar arguments as in Theorem  \ref{classification of schwarzian DE1}, we find that
\begin{itemize}
  \item  all the poles of  $u$  are simple;
  \item $u'(z) = 0$ if and only if $u(z) = \tau_i $  for some $i \in \{1,2,3\}$;
  \item if $u(z) = \tau_1 $, then $z$ is a simple zero of $u'$;
  \item if $u(z) = \tau_j, j=2,3 $, then $z$ is a triple zero of $u'$;
\end{itemize}
It follows that
\begin{equation}
G(z) = \frac{u'^4}{(u-\tau_1)^2(u-\tau_2)^3(u-\tau_3)^3}
\end{equation}
is a nonvanishing entire function, and hence, there exists an entire function $g(z)$ such that $G=e^{g}$.  As the order of $u$ is finite, we have
\begin{eqnarray*}
T(r, e^{g})&=&m(r, e^{g})
\\
&\leq& m\left(r,\frac{u'^2}{(u-\tau_1)^2(u-\tau_2)^2(u-\tau_3)^2}\right)+ m\left(r,\frac{u'}{u-\tau_2}\right)+  m\left(r,\frac{u'}{u-\tau_3}\right)
\\
&\leq& 2m\left(r,\frac{u'}{u-\tau_1}\right)+3m\left(r,\frac{u'}{u-\tau_2}\right)+3m\left(r,\frac{u'}{u-\tau_3}\right)+ O(1)
\\
&=&O(\log r).
\end{eqnarray*}
This indicates that $g= C$ is a constant and hence $u$  satisfies the differential equation
\begin{equation}  \label{subequation-3}
u'^4 = K (u-\tau_1)^2(u-\tau_2)^3(u-\tau_3)^3, \quad  K =  e^C.
\end{equation}
Since the  elliptic curve parametrized by  $u$ and $u'$ has genus one, the general solution of the above equation should be elliptic functions. Then the singularity methods indicate that the general solution to \eqref{subequation-3} can be expressed as
\begin{eqnarray}
   u(z) &=& h+ \frac{1}{2 L} \frac{\wp '\left(z-z_0;g_2,g_3\right)-A_1}{ \wp \left(z-z_0;g_2,g_3\right)-B_1} + \nonumber\\
    &&    \frac{i}{2L} \left(\frac{\wp '\left(z-z_0;g_2,g_3\right)-A_2}{ \wp \left(z-z_0;g_2,g_3\right)-B_2 }-\frac{\wp '\left(z-z_0;g_2,g_3\right)-A_3}{ \wp \left(z-z_0;g_2,g_3\right)-B_3 }\right) \label{solution to Schwarzian DE4}
\end{eqnarray}
where $L^4 = K$, $\wp(z;g_2,g_3)$ is the Weierstrass elliptic function,  $z_0\in\C$ is arbitrary and  $A_j,B_j,g_2,g_3$ are constants  depending on $K$ and $\tau_j,j=1,2,3$. Finally, with $ \tau_1=0, \tau_2=1, \tau_3 =-1 $, substituting \eqref{solution to Schwarzian DE4} into \eqref{ss-eq4} yields that
\begin{eqnarray*}
 && \{\sigma_1, \sigma_2\} = \left\{\frac{i}{2},-\frac{i}{2}\right\},  \quad g_2 = - \frac{c}{36}, \quad B_2 = -B_3 = \frac{L^2}{8}i,
  \\
  &&    c=   \frac{9}{4} L^4, \quad  A_1 = A_2 = A_3 =B_1 = g_3 = h = 0.
\end{eqnarray*}
In this case, the equation \eqref{ss-eq3} reduces to
\begin{eqnarray*}
   S(u,z)^2 = c \frac{(u^2 + 1/4)^2}{u^2 (u^2-1)},
\end{eqnarray*}
and     the solution \eqref{solution to Schwarzian DE4} becomes
\begin{equation*}
u(z) =  -\frac{1}{2 L}\frac{\left(8 \wp \left(z-z_0;g_2,g_3\right)+L^2\right)^2 \wp '\left(z-z_0;g_2,g_3\right)}{ \wp \left(z-z_0;g_2,g_3\right) \left(64 \wp \left(z-z_0;g_2,g_3\right)^2+L^4\right)},
\end{equation*}
where $g_2 = - c/36, g_3 = 0$ and $c=   9 L^4/4$.
This completes the proof.

\begin{remark}
Since elliptic functions are of order $2$,   Theorems \ref{classification of schwarzian DE1}-\ref{classification of schwarzian DE4} indicate that the estimate  on the growth of meromorphic solutions of the equation \eqref{sch-de1} given in Lemma \ref{growth} is sharp.
\end{remark}

\subsection{Examples}
We present some examples to illustrate all the possible configurations of the transcendental meromorphic solutions given in Theorem \ref{classification of schwarzian DE1}.
\begin{example} The Schwarzian differential equation

\begin{eqnarray*}
 S(u,z)=\frac{ 3 \left(25u^4+20 u^3+ 14 u^2+4 u+1\right)}{2 u (u-1) (u+1) \left(3u+1\right) }
\end{eqnarray*}
has the solution
\begin{equation}
  u(z)=\frac{1}{\wp (z-z_0;g_2,g_3)-1},
\end{equation}
where $z_0\in\C$ is arbitrary, $ g_2 = 16$ and $ g_3 = 0.$

\end{example}

\begin{example} The Schwarzian differential equation
\begin{eqnarray*}
  S(u,z)=\frac{ 3 \left(25u^4+20 u^3+ 14 u^2+4 u+1\right)}{ u (u-1) (u+1) \left(3u+1\right) }
\end{eqnarray*}
admits the solution
\begin{equation}
  u(z)=1-\frac{16}{\wp (z-z_0;g_2,g_3) + 12},
\end{equation}
where $z_0\in\C$ is arbitrary, $ g_2 = 64 $ and $ g_3 =0.$

\end{example}

\begin{example} The Schwarzian differential equation

\begin{eqnarray*}
S(u,z) = -\frac{3 \left(25u^4+20 u^3+ 14 u^2+4 u+1\right)}{ u (u-1)  (u+1) (3 u+1)}
\end{eqnarray*}
has the solution
\begin{equation}
  u(z)=-1-\frac{8}{\wp (z-z_0;g_2,g_3) - 8}, 
\end{equation}
where $z_0\in\C$ is arbitrary, $ g_2 = 64$ and $ g_3 = 0.$

\end{example}

\begin{example}

The Schwarzian differential equation
\begin{eqnarray*}
S(u,z) = \frac{3\left(225 u^4+180 u^3+126 u^2+36 u+9\right)}{ u (u-1)  (u+1) (3 u+1)}, 
\end{eqnarray*}
admits the solution
\begin{equation}
  u(z)=-\frac{1}{3}-\frac{16}{\wp (z-z_0;g_2,g_3)-12},
\end{equation}
where $z_0\in\C$ is arbitrary, $ g_2 = 5184$ and $ g_3 = 0.$

\end{example}

\section*{Acknowledgement}
We would like to thank Robert Conte for the helpful discussions.

\section*{Funding}
The first author was supported by the National Natural Science Foundation of China (Grant No. 11671191). The second author was supported by the National Natural Science Foundation of China (Grant
Nos. 11701382 and 11971288).

\vskip 0.3cm
Department of Mathematics\par
Nanjing University\par
Nanjing, China\par
Email: maliao@nju.edu.cn\par

\vskip 0.3cm

Institute for Advanced Study\par
Shenzhen University\par
Shenzhen, China\par
Email: cfwu@szu.edu.cn

\end{document}